\documentclass[11pt]{article}
\usepackage[margin=1in]{geometry}
\usepackage{amsmath,amssymb,amsthm}
\usepackage{mathtools}
\usepackage{enumitem}
\usepackage{hyperref}
\hypersetup{colorlinks=true,linkcolor=blue,citecolor=blue,urlcolor=blue}

\theoremstyle{plain}
\newtheorem{theorem}{Theorem}[section]
\newtheorem{lemma}[theorem]{Lemma}
\newtheorem{proposition}[theorem]{Proposition}
\newtheorem{corollary}[theorem]{Corollary}
\theoremstyle{definition}
\newtheorem{definition}[theorem]{Definition}

\newcommand{\dom}{\operatorname{dom}}
\newcommand{\dist}{\operatorname{dist}}

\usepackage{amsthm}
\newtheorem{question}[theorem]{Question}

\title{Sheaf Subcategories of Quiver Representations}
\author{%
  Eric M. Schmid Jr.\\
  {\small Mathematisch-Naturwissenschaftliche Fakult\"at (Faculty of Science)}\\
  {\small Fachbereich Informatik (Department of Computer Science)}\\
  {\small Eberhard Karls Universit\"at T\"ubingen (University of T\"ubingen)}
  \and
  Fernando Tohm\'e\\
  {\small Universidad Nacional del Sur}\\
  {\small Bah\'ia Blanca, Argentina}
  \and
  William Chin\\
  {\small DePaul University}\\
  {\small Chicago}
}
\date{\today}

\begin{document}
\maketitle

\begin{abstract}
Let $Q$ be a finite quiver without oriented cycles, $C$ its path
category, and $k$ a field. This paper is expository. We assemble in one
place, and in elementary combinatorial form, the dictionary between
Grothendieck topologies on $C$ and subcategories of
$\mathrm{mod}\text{-}kQ$. Every Grothendieck topology on $C$ is rigid, so
the sheaf categories are presheaf categories on full subcategories and are
indexed by subsets $\Sigma$ of the vertices; this is due to Murfet in the
quiver case and, in far greater generality, to Di--Li--Liang. Passing to
$k$-linear coefficients, each topology cuts out a full subcategory
$\mathcal{S}_\Sigma \subseteq \mathrm{mod}\text{-}kQ$. We identify
$\mathcal{S}_\Sigma$ as the perpendicular category, in the sense of
Geigle--Lenzing, of the set of simples off $\Sigma$; it is therefore wide,
and equivalent to $\mathrm{mod}\text{-}kQ_\Sigma$ for an explicit
\emph{reduced quiver} $Q_\Sigma$. We record when $\mathcal{S}_\Sigma$ is a
Serre subcategory (exactly when $\Sigma$ is closed under successors),
observe that $\Sigma \mapsto \mathcal{S}_\Sigma$ embeds the Boolean lattice
$2^{Q_0}$ into the lattice of wide subcategories, and note that in Dynkin
type $A_n$ this captures $2^n$ of the $C_{n+1}$ wide subcategories.
Intrinsically, a wide subcategory $\mathcal{W}$ is a sheaf subcategory
exactly when ${}^{\perp}\mathcal{W}$ is a Serre subcategory and
$\mathcal{W} = ({}^{\perp}\mathcal{W})^{\perp}$. All of
these results are known, or follow readily from known results; the aim is
a concrete self-contained account of the quiver case, with attributions
collected in \S\ref{sec:literature}.
\end{abstract}

\section{Introduction}

The classification of Grothendieck topologies on the path category of a
finite acyclic quiver is settled, and has been settled at several levels
of generality. Grothendieck and Verdier \cite[IV.9]{SGA4} observed that
every Grothendieck topology on a finite Karoubi category is rigid --- a
\emph{subcategory topology} --- so that by the Comparison Lemma
\cite[C2.2.3]{Elephant} the associated sheaf category is again a presheaf
category, on the full subcategory of irreducible objects. Lindenhovius
\cite{Lindenhovius} and Hemelaer \cite{Hemelaer} treated posets;
Murfet \cite{Murfet} treated finite acyclic quivers in an unpublished 2005
note; and Di, Li and Liang \cite{DLL} have recently proved that for a
directed category $\mathcal{C}$ with finite endomorphism monoids, every
Grothendieck topology on $\mathcal{C}^{\mathrm{op}}$ is rigid if and only
if $\mathcal{C}$ is a noetherian EI category. The path category of a
finite acyclic quiver is finite with trivial endomorphism monoids, so it
lies in every one of these classes at once, and the count of topologies is
$2^m$ for $m$ vertices.

We take that classification as given and ask a different question, one
that only makes sense after linearizing. Fix a field $k$ and write $A =
kQ$. A topology $G_\Sigma$, indexed by a subset $\Sigma$ of the vertices,
cuts out a full subcategory
\[
\mathcal{S}_\Sigma \;\subseteq\; \mathrm{mod}\text{-}kQ
\]
of $k$-linear sheaves, which we call a \emph{sheaf subcategory}. Since
$kQ$ is hereditary, $\mathrm{mod}\text{-}kQ$ carries a rich and
well-studied lattice of wide subcategories --- for $Q$ Dynkin these are
counted by generalized Catalan numbers and biject with noncrossing
partitions, clusters and functorially finite torsion classes, by
Ingalls--Thomas \cite{IT}. Our question is:

\begin{question}\label{res:1-1}
Which wide subcategories of $\mathrm{mod}\text{-}kQ$ are sheaf
subcategories?
\end{question}

This is a question about the interaction of two structures: Grothendieck
topologies on one side, the perpendicular calculus of a hereditary algebra
on the other. Both sides are well developed, and the answers below are
consequences of what is already known on each --- we claim no new theorems.
What we offer is the explicit combinatorial dictionary in the quiver case,
which we have not seen written down in one place, together with careful
attributions in \S\ref{sec:literature}.

The answers are as follows. Write $S_v$ for the simple at a vertex $v$ and
$(-)^{\perp}$ for the right perpendicular category in the sense of
Geigle--Lenzing \cite{GL}.

\begin{itemize}
\item (Theorem~\ref{thm:perp}) $\mathcal{S}_\Sigma = \big(\bigoplus_{v
\notin \Sigma} S_v\big)^{\perp}$. Sheaf subcategories are exactly the
perpendicular categories of \emph{sets of simples}. This is the
specialization of \cite[Thm.~4.2]{DLL} to our setting.
\item (Proposition~\ref{prop:wide}) Every $\mathcal{S}_\Sigma$ is a wide
subcategory of $\mathrm{mod}\text{-}kQ$, equivalent to
$\mathrm{mod}\text{-}kQ_\Sigma$ for an explicitly described
\emph{reduced quiver} $Q_\Sigma$, and hence again hereditary --- for
every subset $\Sigma$, with no order condition. Wideness is
\cite[Prop.~1.1]{GL} (using $\operatorname{proj\,dim} S_v \leq 1$, automatic
over the hereditary algebra $kQ$); the equivalence with
$\mathrm{mod}\text{-}eAe$ and the Serre quotient description are
\cite[Prop.~5.1]{GL}, whose functor-category form
\cite[Prop.~5.1$^*$]{GL} is the deletion of vertices.
\item (Proposition~\ref{prop:serre}) $\mathcal{S}_\Sigma$ is a Serre
subcategory if and only if $\Sigma$ is closed under successors --- the
combinatorial form, for path algebras, of the classification of Serre
subcategories of $\mathrm{mod}$ of an Artin algebra in
\cite[Prop.~5.3]{GL} (see also the $\tau$-perpendicular generalization
via \cite{Jasso}). Generic
sheaf subcategories are wide but neither Serre, nor torsion, nor
torsion-free classes.
\item (Proposition~\ref{prop:lattice}) $\Sigma \mapsto G_\Sigma$ is an
anti-isomorphism of the Boolean lattice $2^{Q_0}$ onto the lattice of
topologies, and $\Sigma \mapsto \mathcal{S}_\Sigma$ is an injective
order-preserving map into the lattice of wide subcategories. In Dynkin
type $A_n$ it captures $2^n$ of the $C_{n+1}$ wide subcategories, and we
identify the smallest omitted one. The count $2^n$ appears on the algebra
side already in \cite[Cor.~5.4]{GL}, as the number of surjective
homological epimorphisms out of an Artin algebra with $n$ simples.
\end{itemize}

These answer Question~\ref{res:1-1} by \emph{description}: the sheaf
subcategories are the perpendicular categories of sets of simples. The
\emph{intrinsic} form of the question --- recognizing, from a wide
subcategory alone, whether it is of this shape --- is settled by the
Geigle--Lenzing theory as well: $\mathcal{W}$ is a sheaf subcategory if
and only if $\mathcal{W} = ({}^{\perp}\mathcal{W})^{\perp}$ and the left
perpendicular category ${}^{\perp}\mathcal{W}$ is a Serre subcategory.
See \S\ref{sub:8-4}.

Sections \ref{sec:quivers}--\ref{sec:bijection} recall the classification
of topologies. The results there are Murfet's \cite{Murfet}, and are
subsumed by \cite{DLL}; we include proofs for completeness and because the
quiver case admits a short self-contained treatment, but we claim no
priority for them. Section~\ref{sec:examples} works small examples.
Section~\ref{sec:reps} sets up the linear dictionary and describes the
projectives. Section~\ref{sec:which} contains the results listed above.
Section~\ref{sec:literature} records in detail what is due to whom.

This paper supersedes the preliminary preprint \cite{STC}, which treated
only two extreme topologies and contained errors in the description of
the resulting sheaf categories.

\subsection{Background}

We recall the standard notions the paper builds on; see \cite{MM} for a
textbook treatment of sites and sheaves.

A \textbf{sieve} on an object $c$ of a category $\mathcal C$ is a set $S$
of morphisms with codomain $c$ such that $f \in S$ and $g$ any morphism
with codomain $\dom(f)$ together imply $fg \in S$. Write $t_c$ for the
\textbf{maximal sieve} at $c$, i.e.\ \emph{all} morphisms with codomain
$c$. For a morphism $f : d \to c$ and a sieve $S$ on $c$, the
\textbf{pullback} is
\[
f^*S = \{\, g : x \to d \mid fg \in S \,\}\,,
\]
a sieve on $d$.

A \textbf{Grothendieck topology} $G$ on $\mathcal C$ assigns to each object
$c$ a set $G(c)$ of sieves on $c$ (``covers''), such that:

\begin{itemize}
\item[] \textbf{(Maximality)} $t_c \in G(c)$ for every $c$.
\item[] \textbf{(Stability)} if $S \in G(c)$ and $f : d \to c$, then $f^*S
\in G(d)$.
\item[] \textbf{(Transitivity)} if $S \in G(c)$ and $R$ is a sieve on $c$
such that $f^*R \in G(\dom f)$ for every $f \in S$, then $R \in G(c)$.
\end{itemize}

\section{Quivers, paths, and roots}\label{sec:quivers}

\emph{The material of \S\S\ref{sec:quivers}--\ref{sec:bijection} is due
to Murfet \cite{Murfet}, and is subsumed by the far more general results
of \cite{DLL}; see \S\ref{sec:literature}. We reproduce it, with proofs,
only to keep the paper self-contained, and claim no priority.}

\medskip

A \textbf{quiver} is a nonempty finite directed graph; we call its edges
\textbf{arrows}. A \textbf{composite arrow} is a nonempty sequence of
arrows in which each begins where the previous one ends. A \textbf{cycle}
is a composite arrow beginning and ending at the same vertex. The
\textbf{length} of a composite arrow is the number of arrows it involves.
In a quiver without cycles, a composite arrow visits no vertex twice.

If $Q$ is a quiver, a \textbf{root} is a vertex with no arrows ending at
it; write $r(Q)$ for the set of roots.

Throughout, $Q$ denotes a quiver without cycles and $C$ its \textbf{path
category}: objects are the vertices of $Q$, and for vertices $v,w$ the
morphism set $C(v,w)$ is the empty path (identity) if $v=w$, and otherwise
the set of composite arrows from $v$ to $w$. Write $Q(v,w)$ for the
(direct) arrows $v \to w$. Since $Q$ is finite and acyclic, $C$ has
finitely many morphisms.

\begin{lemma}\label{res:2-1}
If $Q$ is a quiver without cycles, then $r(Q)$ is nonempty.
\end{lemma}

\begin{proof}
For a vertex $v$, let $\ell(v)$ be the maximum length of a composite arrow
beginning at $v$, or $0$ if none exists; this is a well-defined
nonnegative integer since $Q$ is finite and acyclic (paths cannot repeat a
vertex, so their length is bounded by $|Q_0|$). Choose $r$ maximizing
$\ell$. If some arrow $a : u \to r$ existed, then concatenating $a$ with a
longest path $p$ starting at $r$ (possibly $p$ empty, if $\ell(r)=0$)
gives a composite arrow at $u$ of length $\ell(r)+1$, so $\ell(u) \geq
\ell(r)+1 > \ell(r)$ --- contradicting maximality of $\ell(r)$. Hence $r$
has no incoming arrows, i.e.\ $r \in r(Q)$.
\end{proof}

The following is used in Proposition~\ref{res:3-3}'s proof, to induct outward from
the roots.

\begin{definition}[Distance]\label{res:2-2}
Define $\dist : Q_0 \to \mathbb N$ by $\dist(v) = 0$ if $v$ is a root, and
otherwise
\[
\dist(v) = 1 + \max\{\, \dist(u) \mid u \to v \text{ is an arrow of } Q
\,\}\,.
\]
This is well defined: if $a : u \to v$ is an arrow, any path from a root
to $u$ extends via $a$ to a path from a root to $v$ one longer, so
$\dist(v) > \dist(u)$ for every predecessor $u$ of $v$. Since $Q$ is
finite and acyclic (so there is no infinite descending chain of
predecessors), this recursion terminates and defines $\dist$ on every
vertex, with roots exactly the vertices of distance $0$ by Lemma~\ref{res:2-1}.
\end{definition}

\section{Consistent assignments of sieves}\label{sec:consistent}

\begin{definition}\label{res:3-1}
For an arrow $a : v \to w$ and a sieve $S$ at $v$, write $a \cdot S =
\{\, af \mid f \in S \,\}$, a set of morphisms into $w$. A function $J$
assigning to each vertex $v$ a sieve $J(v)$ at $v$ is \textbf{consistent}
if for every vertex $w$, either $J(w) = t_w$, or
\begin{equation}\label{eq:consistent}
J(w) \;=\; \bigcup_{v \in Q_0} \; \bigcup_{a \in Q(v,w)} a \cdot J(v)
\end{equation}
with the convention that the union is $\varnothing$ if $w$ has no incoming
arrows (in particular, at a root, \eqref{eq:consistent} always gives $J(w) = \varnothing$).
\end{definition}

Two standard facts about Grothendieck topologies, needed below.

\begin{lemma}\label{res:3-2}
Let $G$ be a Grothendieck topology on a category with finitely many
morphisms into each object. Then:
\begin{enumerate}[label=(\alph*)]
\item if $S,T \in G(c)$ then $S \cap T \in G(c)$;
\item if $S \in G(c)$ and $T$ is a sieve on $c$ with $T \supseteq S$, then
$T \in G(c)$.
\end{enumerate}
Consequently $J_G(c) := \bigcap_{S \in G(c)} S$ lies in $G(c)$ and is its
minimum element, and $G(c) = \{\, T \text{ sieve on } c \mid T \supseteq
J_G(c) \,\}$.
\end{lemma}

\begin{proof}
(b) For $f \in S \subseteq T$, since $T$ is a sieve, $f^*T =
t_{\dom f}$ (everything precomposes into $T$ once $f$ does); and
$t_{\dom f} \in G(\dom f)$ always. So $f^*T \in G(\dom f)$
for every $f \in S$, and transitivity (applied to the cover $S$ and
candidate sieve $T$) gives $T \in G(c)$.

(a) For $f \in T$, as above $f^*T = t_{\dom f}$, so $f^*(S \cap T) = f^*S
\cap t_{\dom f} = f^*S$, which lies in $G(\dom f)$ by stability (since $S
\in G(c)$). So $f^*(S\cap T) \in G(\dom f)$ for every $f \in T$, and
transitivity (applied to the cover $T$ and candidate $S \cap T$) gives
$S \cap T \in G(c)$.

For the consequence: finitely many covers at $c$ means iterating (a)
produces $J_G(c) \in G(c)$; it is a lower bound for $G(c)$ by
construction, so it is the minimum; every $T \supseteq J_G(c)$ is then a
cover by (b), and every cover trivially contains $J_G(c)$.
\end{proof}

\begin{proposition}\label{res:3-3}
If $Q$ has $m$ vertices, there are exactly $2^m$ consistent assignments of
sieves.
\end{proposition}

\begin{proof}
We biject consistent assignments with functions $g : Q_0 \to \{0,1\}$.

\emph{From $g$ to $J_g$.} For a root $r$, the only sieves at $r$ are
$\{1_r\} = t_r$ and $\varnothing$ (since $r$ has no incoming arrows, the
only morphism into $r$ at all is $1_r$). Set $J_g(r) = t_r$ if $g(r)=0$,
and $J_g(r) = \varnothing$ if $g(r)=1$ --- note this matches \eqref{eq:consistent},
which at a root always evaluates to $\varnothing$. Proceeding by induction
on $\dist$ (well founded by Definition~\ref{res:2-2}, and using that every
predecessor of $w$ has strictly smaller distance than $w$), define, once
$J_g$ is known on all vertices of distance $< n$,
\begin{equation}\label{eq:extend}
J_g(w) = \begin{cases} t_w & g(w) = 0 \\[4pt] \displaystyle
\bigcup_{v,\, a \in Q(v,w)} a \cdot J_g(v) & g(w) = 1 \end{cases}
\end{equation}
for $w$ of distance $n$. By construction $J_g$ satisfies \eqref{eq:consistent} at every
vertex, i.e.\ is consistent.

\emph{From $J$ to $g_J$.} Given consistent $J$, set $g_J(w) = 0$ if $J(w)
= t_w$, else $g_J(w)=1$.

\emph{These are mutually inverse.} First, $J_{g_J} = J$: at a root $r$
this is immediate (only two sieves exist there, and the definitions match
up by cases). Inductively, suppose $J_{g_J}(v) = J(v)$ for all $v$ of
distance $< n$, and let $w$ have distance $n$. If $J(w) = t_w$ then
$g_J(w)=0$ and $J_{g_J}(w) = t_w = J(w)$ by \eqref{eq:extend}. If $J(w) \neq t_w$ then
$g_J(w)=1$, and consistency of $J$ gives $J(w) = \bigcup a \cdot J(v)$
(sum over predecessors $v$, all of distance $<n$), which by the inductive
hypothesis equals $\bigcup a \cdot J_{g_J}(v) = J_{g_J}(w)$ by \eqref{eq:extend}. This
closes the induction.

Second, $g_{J_g} = g$: note that if $g(w)=1$, then $J_g(w) = \bigcup a
\cdot J_g(v)$ never contains the identity $1_w$ --- because $a \cdot
J_g(v)$ consists of composite arrows of length $\geq 1$ ending at $w$, and
since $Q$ is acyclic no such arrow can equal $1_w$ (that would be a
cycle). Since $t_w$ \emph{does} contain $1_w$, this shows $J_g(w) \ne
t_w$ whenever $g(w)=1$, i.e.\ $J_g(w)$ is a proper sieve. Hence
\[
g(w) = 0 \iff J_g(w) = t_w \iff g_{J_g}(w) = 0,
\]
the middle equivalence by definition of $J_g$ and the fact just shown, the
last by definition of $g_{J_g}$. So $g = g_{J_g}$.

This is a bijection between $\{0,1\}^{Q_0}$ (size $2^m$) and consistent
assignments.
\end{proof}

\section{The bijection with Grothendieck topologies}\label{sec:bijection}

The following lemma isolates the combinatorial heart of the
transitivity axiom; it depends only on consistency of $J$, and is
used in the proof of Theorem~\ref{res:4-1} below.

\begin{lemma}\label{lem:transitivity}
Let $J$ be consistent and $w$ a vertex with $J(w) \neq t_w$ and $J(w) \ne
\varnothing$. Let $T \supseteq J(w)$ be a sieve, and $S$ a sieve at $w$
with $h^*S \supseteq J(\dom h)$ for every $h \in T$. Then $J(w) \subseteq
S$.
\end{lemma}

\begin{proof}
Fix $f \in J(w)$; we find a prefix $h$ of $f$, lying in $T$, whose
``remaining suffix'' lies in $J(\dom h)$, and conclude from there.

Since $1_w \notin J(w)$ (identity morphisms never lie in the union \eqref{eq:consistent}, by
the acyclicity argument of Proposition~\ref{res:3-3}), $f$ is a genuine nonempty
path. By consistency, $J(w) = \bigcup a \cdot J(v)$, so $f$ factors
uniquely as $f = a_1 g_1$ where $a_1$ is the final arrow of $f$ and $g_1
\in J(v_1)$, $v_1 = \dom(a_1)$.

\begin{itemize}
\item If $J(v_1) = t_{v_1}$: then $1_{v_1} \in J(v_1)$, so $a_1 = a_1
\cdot 1_{v_1} \in a_1 \cdot J(v_1) \subseteq J(w) \subseteq T$. Take $h =
a_1$; then $f = h g_1$ with $h \in T$ and $g_1 \in J(v_1) = J(\dom h)$.
Stop.
\item Otherwise $J(v_1) \neq t_{v_1}$; since $g_1 \in J(v_1)$ witnesses
$J(v_1) \neq \varnothing$, consistency gives $J(v_1) = \bigcup a \cdot
J(v)$ again, so $g_1 = a_2 g_2$ with $a_2$ the final arrow of $g_1$ and
$g_2 \in J(v_2)$, $v_2 = \dom(a_2)$. Repeat the test on $J(v_2)$, now with
the longer prefix $a_1a_2$.
\end{itemize}

Each iteration strictly shortens the remaining suffix, so after at most
$\ell(f)$ steps the process must stop. It cannot fail to stop: if it were
still running when the remaining suffix reached the identity $1_{v_k}$ at
$v_k = \dom(f)$, that identity would have to lie in $J(v_k)$ (it is the
witnessing element at that stage) --- but $1_{v_k} \in J(v_k)$ forces
$J(v_k) = t_{v_k}$ exactly (a sieve containing an identity is the maximal
sieve), which \emph{is} the stopping condition. So the process necessarily
halts at some $h = a_1 \cdots a_k \in T$ with remaining suffix $g_k \in
J(\dom h)$ and $f = h g_k$.

Now $h \in T$ gives $h^*S \supseteq J(\dom h) \ni g_k$, so $g_k \in h^*S$,
i.e.\ $h g_k = f \in S$.
\end{proof}

\begin{theorem}\label{res:4-1}
There is a bijection between Grothendieck topologies on $C$ and consistent
assignments of sieves.
\end{theorem}

\begin{proof}
\textbf{($G \mapsto J_G$.)} Given a topology $G$, let $J_G(v) =
\bigcap_{S \in G(v)} S$ as in Lemma~\ref{res:3-2}, so $G(v) = \{\, T \mid T
\supseteq J_G(v)\,\}$. We check $J_G$ is consistent.

At a root $r$: the only sieves are $t_r,\varnothing$, so either $G(r) =
\{t_r\}$ (giving $J_G(r)=t_r$) or $G(r) = \{t_r,\varnothing\}$ (giving
$J_G(r)=\varnothing$) --- either way \eqref{eq:consistent} holds trivially at $r$.

At a non-root $w$: if $G(w) = \{t_w\}$ then $J_G(w)=t_w$ and we are done.
Otherwise, let $T = \bigcup_{v,\,a\in Q(v,w)} a \cdot J_G(v)$; we must show
$J_G(w) = T$.

\emph{$J_G(w) \supseteq T$:} fix an arrow $a : v \to w$. By stability,
$a^* J_G(w) \in G(v)$; since $J_G(v)$ is the \emph{minimum} cover at $v$,
this forces $a^*J_G(w) \supseteq J_G(v)$, i.e.\ for every $g \in J_G(v)$,
$ag \in J_G(w)$. This holds for every arrow into $w$, so $T \subseteq
J_G(w)$.

\emph{$T \supseteq J_G(w)$, hence equality:} since $T \subseteq J_G(w)$
already, it suffices (by minimality of $J_G(w)$ among covers) to show $T$
itself is a cover, i.e.\ $T \in G(w)$. By transitivity applied to the
cover $J_G(w)$, it is enough to show $f^*T \in G(\dom f)$ for every $f \in
J_G(w)$.

Since $G(w) \neq \{t_w\}$, $J_G(w)$ is proper, so such $f$ is a nonempty
path ending at $w$ via a specific last arrow $b : v \to w$; write $f = b
f'$. For any \emph{other} arrow $a : q \to w$ ($a\ne b$), no morphism $h$
can satisfy $fh = ag$ for $g \in J_G(q)$ (the two sides end in different
final arrows, $b$ and $a$, and $C$ is freely generated by $Q$), so
$f^*(a \cdot J_G(q)) = \varnothing$. Hence only the $(v,b)$ term of $T$
survives the pullback:
\[
f^*T = f^*(b \cdot J_G(v)) = (f')^* J_G(v)
\]
(cancel the common initial factor $b$). Since $J_G(v)$ is a cover at $v$
and $f' : \dom f \to v$ is a morphism, stability gives $(f')^*J_G(v) \in
G(\dom f)$, i.e.\ $f^*T \in G(\dom f)$, as required. So $T \in G(w)$,
giving $J_G(w) = T$ and consistency of $J_G$.

\textbf{($J \mapsto G_J$.)} Given consistent $J$, set $G_J(v) = \{\, T
\mid T \supseteq J(v)\,\}$. Maximality is immediate ($t_v \supseteq J(v)$
always).

\emph{Stability.} For an arrow $a : v \to w$, consistency gives $a \cdot
J(v) \subseteq J(w)$ (either because $J(w)=t_w$, trivially, or because
$a\cdot J(v)$ is literally one term of the union \eqref{eq:consistent} defining $J(w)$);
unwinding definitions this says $a^*J(w) \supseteq J(v)$. Composing this
fact along a general morphism $f = a_n \cdots a_1$ gives $f^*J(w)
\supseteq J(v)$ for \emph{any} morphism $f : v \to w$ of $C$ (pullback is
monotone and functorial). Hence if $S \in G_J(w)$, i.e.\ $S \supseteq
J(w)$, then $f^*S \supseteq f^*J(w) \supseteq J(v)$, i.e.\ $f^*S \in
G_J(v)$.

\emph{Transitivity.} Suppose $T \supseteq J(w)$ and $S$ is a sieve at $w$
with $h^*S \in G_J(\dom h)$ (i.e.\ $h^*S \supseteq J(\dom h)$) for every $h
\in T$; we must show $S \supseteq J(w)$, i.e.\ $S \in G_J(w)$. If $J(w)
\in \{t_w, \varnothing\}$ this is immediate ($J(w)=t_w$ forces $T=t_w \ni
1_w$, giving $S = 1_w^*S \supseteq J(w)$ directly; $J(w)=\varnothing$
needs nothing). Otherwise it is exactly Lemma~\ref{lem:transitivity}. So
$G_J$ is a Grothendieck topology.

\textbf{(Mutual inverses.)} $J_{G_J} = J$: by construction $G_J(v) = \{T
\mid T \supseteq J(v)\}$ is the up-set of $J(v)$, whose intersection is
$J(v)$ itself. $G_{J_G} = G$: this is the identity $G(v) = \{T \mid T
\supseteq J_G(v)\}$ established in Lemma~\ref{res:3-2}, applied to $G$.
\end{proof}

\begin{corollary}\label{res:4-3}
Let $Q$ be a quiver without cycles with $m$ vertices, and $C$ its path
category. There are exactly $2^m$ Grothendieck topologies on $C$, in
bijection with subsets of vertices via
\[
G \longmapsto \{\, v \in Q_0 \mid G(v) = \{t_v\} \,\}.
\]
\end{corollary}

\begin{proof}
Immediate from Proposition~\ref{res:3-3} and Theorem~\ref{res:4-1}, tracking through that
$g(v) = 0 \iff J_G(v) = t_v \iff G(v) = \{t_v\}$.
\end{proof}

\section{Worked examples}\label{sec:examples}

\subsection{A single arrow}\label{sub:5-1}

Take $Q$ with two vertices $1,2$ and a single arrow $a : 1 \to 2$. Here
$m=2$, so Corollary~\ref{res:4-3} predicts $2^2 = 4$ topologies. Vertex $1$ is the
unique root.

Sieves at $1$: $\varnothing$ or $t_1 = \{1_1\}$ --- two choices. Sieves at
$2$: a sieve containing $1_2$ must (by the sieve axiom, since $1_2 a = a$)
also contain $a$, so the only sieves are $\varnothing$, $\{a\}$, and $t_2
= \{a,1_2\}$ --- three choices, of which \eqref{eq:consistent} can only ever produce
$\varnothing$ (if $J(1)=\varnothing$) or $\{a\}$ (if $J(1) = t_1$).

The four consistent assignments are exactly $(J(1),J(2)) \in
\{(t_1,t_2),\ (t_1,\{a\}),\ (\varnothing,t_2),\ (\varnothing,\varnothing)\}$,
corresponding to the four topologies
\[
\begin{array}{llll}
G(1)=\{t_1\}, & G(2) = \{t_2\} & \text{(chaotic)}\\[4pt]
G(1)=\{t_1\}, & G(2) = \{\{a\}, t_2\}\\[4pt]
G(1)=\{\varnothing,t_1\}, & G(2) = \{t_2\}\\[4pt]
G(1)=\{\varnothing,t_1\}, & G(2) = \{\varnothing,\{a\},t_2\} & \text{(all sieves cover)}
\end{array}
\]
One can check directly (stability along $a$, transitivity) that these are
exactly the four valid choices --- in particular the combination $G(1)
\ni \varnothing$ together with $G(2) = \{\{a\},t_2\}$ (missing
$\varnothing$) \emph{fails} transitivity, correctly excluded by the
theorem.

\subsection{The Kronecker quiver}\label{sub:5-2}

Take $Q$ with two vertices $1,2$ and \emph{two} parallel arrows $a,b : 1
\to 2$ (the Kronecker quiver). This is still acyclic, and still has $m=2$
vertices, so Corollary~\ref{res:4-3} again predicts $2^2 = 4$ topologies --- the
count depends only on the number of vertices, not on how many arrows
connect them. It's instructive to see why, since the local sieve lattice
here is strictly bigger than in \S\ref{sub:5-1}.

Sieves at $1$: as before, $\varnothing$ or $t_1$ --- two choices. Sieves
at $2$: a subset of $\{a,b,1_2\}$ is a sieve iff it either omits $1_2$
entirely (any of the $2^2=4$ subsets of $\{a,b\}$ works, since $a,b$ both
have domain $1$ and precomposing with $1_1$ changes nothing), or contains
$1_2$, in which case closure forces it to contain $a$ and $b$ too, making
it all of $t_2$. So there are $5$ sieves at $2$: $\varnothing$, $\{a\}$,
$\{b\}$, $\{a,b\}$, and $t_2$.

Formula (1) at vertex $2$ reads $J(2) = a\cdot J(1) \cup b \cdot J(1)$.
Since $J(1) \in \{\varnothing, t_1\}$, this union is always either
$\varnothing$ (if $J(1)=\varnothing$) or $\{a,b\}$ (if $J(1)=t_1$) ---
\emph{never} $\{a\}$ or $\{b\}$ alone. So, exactly as in \S\ref{sub:5-1}, $J(2)$
has only two possible values once $J(1)$ is fixed: $t_2$, or the formula
value. This gives $2 \times 2 = 4$ consistent assignments, matching
Corollary~\ref{res:4-3}:
\[
(J(1),J(2)) \in \{(t_1,t_2),\ (t_1,\{a,b\}),\ (\varnothing,t_2),\
(\varnothing,\varnothing)\}.
\]
Two of the five sieves at vertex $2$ --- namely $\{a\}$ and $\{b\}$ ---
never occur as $J_G(2)$ for \emph{any} Grothendieck topology $G$ on $C$:
a topology here cannot cover $2$ by one of the two arrows alone without
also covering it by the other. This is a general feature of the
correspondence, not a special coincidence of this example: \eqref{eq:consistent}
only ever produces the \emph{full} pushforward $\bigcup_a a\cdot J(v)$, so
any topology on a quiver's path category treats all arrows into a vertex
symmetrically at the level of what can be a minimal cover.

\subsection{Type $A_3$, equioriented}\label{sub:5-3}

Take $Q$ the linearly oriented quiver of type $A_3$: vertices $1,2,3$
with arrows $a : 1\to 2$ and $b : 2 \to 3$. Here $m=3$, so Corollary~\ref{res:4-3}
predicts $2^3=8$ topologies. Vertex $1$ is the unique root;
$\dist(2)=1$, $\dist(3)=2$.

Sieves at $1$: $\varnothing$ or $t_1$ (two, as always at a root). Sieves
at $2$: $\varnothing$, $\{a\}$, or $t_2$ (three, exactly as in \S\ref{sub:5-1},
since $2$ has a single incoming arrow). Sieves at $3$: the morphisms into
$3$ are $\{ba, b, 1_3\}$ (the composite $ba:1\to3$, the arrow $b:2\to3$,
and the identity), and one checks directly that the sieves form a
\emph{chain}
\[
\varnothing \;\subset\; \{ba\} \;\subset\; \{ba,b\} \;\subset\; t_3 =
\{ba,b,1_3\},
\]
since containing $b$ forces containing $ba$ (precompose $b$ with $a$),
and containing $1_3$ forces containing both. So there are $4$ sieves at
vertex $3$. They form a chain here, rather than the lattice seen at the
Kronecker quiver's vertex $2$ in \S\ref{sub:5-2}, because vertex $3$ has
two non-identity morphisms into it, $b$ and $ba$, but only one direct
arrow $b \in Q(2,3)$; it is parallel arrows that produce incomparable
sieves.

Working outward: $J(2) \in \{t_2, a\cdot J(1)\}$, where $a \cdot J(1)$ is
$\{a\}$ if $J(1)=t_1$ and $\varnothing$ if $J(1)=\varnothing$ --- two
choices, as in \S\ref{sub:5-1}. Then $J(3) \in \{t_3,\ b \cdot J(2)\}$, where
\[
b \cdot J(2) = \begin{cases}
\{ba, b\} & J(2) = t_2 \\
\{ba\} & J(2) = \{a\} \\
\varnothing & J(2) = \varnothing
\end{cases}
\]
again two choices once $J(2)$ is fixed. In total $2\times 2\times 2 = 8$
consistent assignments, matching Corollary~\ref{res:4-3}:
\[
\begin{array}{llll}
1.&(t_1,\,t_2,\,t_3) & 5.&(\varnothing,\,t_2,\,t_3)\\
2.&(t_1,\,t_2,\,\{ba,b\}) & 6.&(\varnothing,\,t_2,\,\{ba,b\})\\
3.&(t_1,\,\{a\},\,t_3) & 7.&(\varnothing,\,\varnothing,\,t_3)\\
4.&(t_1,\,\{a\},\,\{ba\}) & 8.&(\varnothing,\,\varnothing,\,\varnothing)
\end{array}
\]
Unlike the Kronecker quiver, here \emph{every} sieve in each vertex's
lattice does arise as some $J_G$: the four sieves at vertex $3$ listed
above are exactly the four values $J(3)$ takes across these eight rows
($t_3$, $\{ba,b\}$, $\{ba\}$, $\varnothing$, each occurring twice). The
difference from \S\ref{sub:5-2} is that $A_3$ has no vertex with two \emph{parallel}
arrows feeding into it, so there is no pair of ``partial'' sieves like
$\{a\},\{b\}$ that \eqref{eq:consistent} is forced to skip.

\medskip
\noindent\emph{Remark.} The other two orientations of the $A_3$ diagram
($1 \to 2 \leftarrow 3$ and $1 \leftarrow 2 \to 3$) still have $m=3$
vertices, so Corollary~\ref{res:4-3} gives $8$ topologies for those quivers as well
--- even though they have two roots instead of one, and the local sieve
counts at each vertex differ from the equioriented case. The count
$2^m$ never sees the orientation, only the number of vertices.

\section{Quiver representations and the categories of projectives}\label{sec:reps}

The classification of \S\ref{sec:bijection} is a statement about
categories of sets. In this section we linearize it, fixing a field $k$,
and describe the resulting subcategories of $\mathrm{mod}\text{-}kQ$ and
their projectives. Theorem~\ref{thm:comparison} below is the
specialization to path categories of the Comparison Lemma
\cite[C2.2.3]{Elephant}, and of \cite[Thm.~5.18]{DLL}; its $k$-linear
form is \cite[Prop.~5.1$^*$]{GL}. We give the explicit combinatorial
version, which is what \S\ref{sec:which} needs.

\subsection{Conventions}\label{sub:6-1}

Murfet's sieves consist of morphisms \emph{into} an object, so the
natural home for the classification is the presheaf category
$\mathrm{PSh}(C) = \mathrm{Fun}(C^{\mathrm{op}}, \mathbf{Set})$: an object
$F$ assigns a set to each vertex and, to each arrow $a : v \to w$, a map
$F(a) : F(w) \to F(v)$ running \emph{against} the arrow. Linearizing over
a field $k$ gives $\mathrm{Fun}(C^{\mathrm{op}},\mathrm{Vect}_k) \simeq
\mathrm{Mod}\text{-}kQ$, representations of $Q^{\mathrm{op}}$ in the
covariant convention. Readers who prefer maps following the arrows should
apply everything below to $Q^{\mathrm{op}}$; nothing else changes. From
\S\ref{sub:6-6} onwards we work with finite-dimensional representations
and write $\mathrm{mod}\text{-}kQ$; since $Q$ is finite and acyclic,
$kQ$ is finite-dimensional and every statement below has an evident
analogue for $\mathrm{Mod}\text{-}kQ$, but Krull--Schmidt is used in
Corollary~\ref{res:6-9} and we prefer to fix the finite-dimensional
convention once.

For a path $f : u \to w$ in $Q$ write $V^{\circ}(f)$ for the set of
vertices \emph{visited by $f$ other than its target $w$} (so $u \in
V^{\circ}(f)$, and $V^{\circ}(1_w) = \varnothing$). Fix throughout a
subset $\Sigma \subseteq Q_0$ and let $G_\Sigma$ be the corresponding
topology from Corollary~\ref{res:4-3}, so that $J(v) = t_v$ exactly for $v \in
\Sigma$. (We write $\Sigma$ rather than $S$ to avoid a clash with the
letters used for sieves in \S\S\ref{sec:consistent}--\ref{sec:bijection}.)

We use the standard fact that representables $y(w) = C(-,w)$ are
projective in $\mathrm{PSh}(C)$: by Yoneda $\mathrm{Hom}(y(w), F) \cong
F(w)$, and evaluation at $w$ preserves epimorphisms because epis in
$\mathrm{PSh}(C)$ are computed objectwise. Linearizing this \emph{is} the
classical statement that $P(w) = k\,C(-,w)$ are the indecomposable
projectives of $\mathrm{Mod}\text{-}kQ$, with $kQ = \bigoplus_w P(w)$.

\subsection{The covering sieves, explicitly}\label{sub:6-2}

\begin{lemma}\label{res:6-1}
Let $w \in Q_0$. If $w \in \Sigma$ then $J(w) = t_w$. If $w \notin \Sigma$
then
\[
J(w) = \{\, f \text{ a path into } w \;\mid\; V^{\circ}(f) \cap \Sigma
\neq \varnothing \,\},
\]
i.e.\ $J(w)$ consists of exactly those paths into $w$ that visit a vertex
of $\Sigma$ strictly before arriving at $w$.
\end{lemma}

\begin{proof}
Induction on $\dist(w)$. If $w$ is a root and $w \notin \Sigma$, then
$J(w) = \varnothing$; correspondingly the only path into $w$ is $1_w$,
with $V^{\circ}(1_w) = \varnothing$, so the right-hand side is empty too.

For the inductive step let $w \notin \Sigma$, so by consistency $J(w) =
\bigcup_{a : v \to w} a \cdot J(v)$. A path $f \in J(w)$ is exactly one of
the form $f = ag$ with $a : v \to w$ its final arrow and $g \in J(v)$, and
$V^{\circ}(f) = V^{\circ}(g) \cup \{v\}$. If $v \in \Sigma$ then $J(v) =
t_v$ imposes no condition on $g$, and indeed $V^{\circ}(f) \ni v$ meets
$\Sigma$. If $v \notin \Sigma$ then by the inductive hypothesis $g \in
J(v)$ iff $V^{\circ}(g) \cap \Sigma \neq \varnothing$, which since $v
\notin \Sigma$ is the same condition as $V^{\circ}(f) \cap \Sigma \neq
\varnothing$.
\end{proof}

\begin{lemma}\label{res:6-2}
A presheaf $F$ is a $G_\Sigma$-sheaf if and only if $F$ satisfies the
sheaf condition with respect to the minimal cover $J(v)$ at every vertex
$v$ simultaneously.
\end{lemma}

\begin{proof}
Necessity is trivial. For sufficiency, assume $F$ satisfies the sheaf
condition for $J(v)$ at every $v$, and let $T \supseteq J(v)$ be any
cover. \emph{Separatedness}: if $x,y \in F(v)$ satisfy $x \cdot f = y
\cdot f$ for all $f \in T$, then in particular for all $f \in J(v)$,
whence $x = y$ by the uniqueness half of the $J(v)$-condition. This holds
at every vertex.

\emph{Existence}: let $(s_f)_{f \in T}$ be a matching family. Restricting
to $J(v) \subseteq T$ and amalgamating gives a unique $x \in F(v)$ with $x
\cdot f = s_f$ for $f \in J(v)$. For arbitrary $f \in T$ consider the
cover $f^*J(v)$ of $\dom f$ (Stability). For $h \in f^*J(v)$ we have $fh
\in J(v)$, so $(x\cdot f)\cdot h = x \cdot (fh) = s_{fh} = s_f \cdot h$.
Thus $x \cdot f$ and $s_f$ agree on the cover $f^*J(v)$, so $x\cdot f =
s_f$ by separatedness at $\dom f$.
\end{proof}

\subsection{The reduced quiver}\label{sub:6-3}

\begin{definition}\label{res:6-3}
The \textbf{reduced quiver} $Q_\Sigma$ has vertex set $\Sigma$, and an
arrow $u \to w$ for each nonempty path $u \to w$ in $Q$ all of whose
\emph{intermediate} vertices (those other than $u$ and $w$) lie outside
$\Sigma$. For $w \notin \Sigma$, let $\mathrm{Red}_\Sigma(w)$ denote the
set of paths $\pi : u \to w$ in $Q$ with $u \in \Sigma$ and all vertices
strictly between $u$ and $w$ outside $\Sigma$; write $s(\pi) = u$. Since
$Q$ is finite and acyclic, so is $Q_\Sigma$.
\end{definition}

\begin{lemma}\label{res:6-4}
Let $C|_\Sigma$ be the full subcategory of $C$ on the objects $\Sigma$.
Then $C|_\Sigma$ is the path category of $Q_\Sigma$.
\end{lemma}

\begin{proof}
Morphisms of $C|_\Sigma$ are the paths of $Q$ between $\Sigma$-vertices.
Given such a path $f : u \to w$, let $u = x_0, x_1, \dots, x_r = w$ be the
$\Sigma$-vertices visited by $f$, in order. Splitting $f$ at these
vertices expresses it as a composite of $r$ arrows of $Q_\Sigma$, and this
factorization is unique: any factorization of $f$ in $C|_\Sigma$ splits
$f$ at a subset of its $\Sigma$-vertices, and a factor is an arrow of
$Q_\Sigma$ precisely when it contains no $\Sigma$-vertex in its interior,
which forces the splitting set to be all of $\{x_0,\dots,x_r\}$. Hence
$C|_\Sigma$ is free on $Q_\Sigma$.
\end{proof}

\begin{proposition}\label{res:6-5}
For a presheaf $F$ on $C$ the following are equivalent.
\begin{enumerate}[label=(\roman*)]
\item $F$ is a $G_\Sigma$-sheaf;
\item \emph{(unrolled form)} for every $w \notin \Sigma$ the canonical map
\[
F(w) \longrightarrow \prod_{\pi \in \mathrm{Red}_\Sigma(w)} F(s(\pi)),
\qquad x \longmapsto (x \cdot \pi)_\pi,
\]
is a bijection;
\item \emph{(one-step form)} for every $w \notin \Sigma$ the canonical map
\[
\theta^F_w : F(w) \longrightarrow \prod_{a\,:\,v \to w} F(v), \qquad x
\longmapsto (x\cdot a)_a,
\]
assembled from the structure maps along the arrows into $w$, is a
bijection.
\end{enumerate}
In particular, if no vertex of $\Sigma$ admits a path to $w$, then $F(w)$
is a singleton.
\end{proposition}

\begin{proof}
By Lemma~\ref{res:6-2} it suffices to analyse the condition $F(w) \cong
\mathrm{Match}(J(w),F)$ for $w \notin \Sigma$ (at $w \in \Sigma$ the
minimal cover is $t_w$ and the condition is vacuous).

Because $C$ is a \emph{free} path category, every $f \in J(w)$ factors
uniquely as $f = ag$ with $a$ its final arrow, and $fh = a(gh)$ has the
same final arrow. Hence a matching family on $J(w) = \bigcup_{a : v \to w}
a \cdot J(v)$ splits, with no cross-conditions, into one matching family
on $J(v)$ for each arrow $a : v \to w$:
\[
\mathrm{Match}(J(w),F) \;\cong\; \prod_{a : v \to w}
\mathrm{Match}(J(v),F).
\]
Now iterate. Each factor with $v \in \Sigma$ contributes
$\mathrm{Match}(t_v,F) = F(v)$ and terminates; each factor with $v \notin
\Sigma$ expands again. The iteration terminates because $Q$ is finite and
acyclic, and a branch that reaches a root outside $\Sigma$ contributes
$\mathrm{Match}(\varnothing,F) = \{*\}$, i.e.\ nothing. The surviving
branches are exactly the paths of $\mathrm{Red}_\Sigma(w)$, giving the
displayed product. (If there are none, the empty product is a singleton.)
This proves (i) $\Leftrightarrow$ (ii).

(ii) $\Rightarrow$ (iii). Assume (ii) and induct on $\dist(w)$ for $w
\notin \Sigma$. Every arrow $a : v \to w$ has $\dist(v) < \dist(w)$.
Splitting $\mathrm{Red}_\Sigma(w)$ according to the final arrow $a$ of
each $\pi$, the paths with final arrow $a : v\to w$ are: the single path
$a$ itself when $v \in \Sigma$; and the paths $a\pi'$ with $\pi' \in
\mathrm{Red}_\Sigma(v)$ when $v \notin \Sigma$. Hence
\[
\prod_{\pi \in \mathrm{Red}_\Sigma(w)} F(s(\pi)) \;\cong\;
\prod_{a : v \to w,\ v \in \Sigma} F(v) \;\times\!\!
\prod_{a : v \to w,\ v \notin \Sigma} \ \prod_{\pi' \in
\mathrm{Red}_\Sigma(v)} F(s(\pi')),
\]
and by the inductive hypothesis each inner product is $F(v)$. The
resulting bijection is compatible with the canonical maps, so (ii) at $w$
becomes exactly the assertion that $\theta^F_w$ is a bijection.

(iii) $\Rightarrow$ (ii). The same computation read in the other
direction: substituting $F(v) \cong \prod_{\pi' \in
\mathrm{Red}_\Sigma(v)} F(s(\pi'))$ for each predecessor $v \notin
\Sigma$ (inductive hypothesis) into $\theta^F_w$ recovers the product
over $\mathrm{Red}_\Sigma(w)$.
\end{proof}

\subsection{The comparison theorem}\label{sub:6-4}

\begin{theorem}\label{thm:comparison}\label{res:6-6}
Restriction along $C|_\Sigma \hookrightarrow C$ is an equivalence
\[
\mathrm{Sh}(C, G_\Sigma) \;\xrightarrow{\ \sim\ }\;
\mathrm{PSh}(C|_\Sigma) \;=\; \mathrm{PSh}(Q_\Sigma).
\]
\end{theorem}

\begin{proof}
We verify the hypotheses of Grothendieck's Comparison Lemma for the full
subcategory $C|_\Sigma \subseteq C$.

\emph{$C|_\Sigma$ is $G_\Sigma$-dense.} Let $R_w$ be the sieve at $w$
generated by all morphisms into $w$ with domain in $\Sigma$, i.e.\ the set
of paths into $w$ that factor through a $\Sigma$-vertex. If $w \in \Sigma$
then $1_w \in R_w$, so $R_w = t_w$ is covering. If $w \notin \Sigma$, then
a path into $w$ factors through a $\Sigma$-vertex precisely when
$V^{\circ}$ of it meets $\Sigma$, so $R_w = J(w)$ by Lemma~\ref{res:6-1} --- again
covering.

\emph{The induced topology on $C|_\Sigma$ is trivial.} A sieve $T$ at $u
\in \Sigma$ covers in the induced topology iff the sieve it generates in
$C$ contains $J(u) = t_u$, i.e.\ iff $1_u$ factors as $fh$ with $f \in T$.
In a path category on an acyclic quiver the only split epimorphisms are
identities, so this forces $f = 1_u \in T$, i.e.\ $T = t_u$. Hence only
maximal sieves cover, and $\mathrm{Sh}(C|_\Sigma, \text{induced}) =
\mathrm{PSh}(C|_\Sigma)$.

The Comparison Lemma now gives the equivalence, and Lemma~\ref{res:6-4} identifies
$C|_\Sigma$ as the path category of $Q_\Sigma$.
\end{proof}

Proposition~\ref{res:6-5} makes the inverse equivalence explicit: given $F_0 \in
\mathrm{PSh}(Q_\Sigma)$, extend it to $C$ by $F(w) = \prod_{\pi \in
\mathrm{Red}_\Sigma(w)} F_0(s(\pi))$ for $w \notin \Sigma$.

\subsection{Projectives}\label{sub:6-5}

\begin{corollary}\label{res:6-7}
Let $a : \mathrm{PSh}(C) \to \mathrm{Sh}(C,G_\Sigma)$ be sheafification.
Then:
\begin{enumerate}[label=(\alph*)]
\item The objects $a\,y(u)$ for $u \in \Sigma$ are projective and
generate $\mathrm{Sh}(C,G_\Sigma)$; under the equivalence of
Theorem~\ref{res:6-6} they correspond to the representables $y_{Q_\Sigma}(u)$.
\item For $w \notin \Sigma$ there is an isomorphism
\[
a\,y(w) \;\cong\; \coprod_{\pi \in \mathrm{Red}_\Sigma(w)} a\,y(s(\pi)),
\]
so $a\,y(w)$ is again projective, but contributes no new generator ---
and is the initial object when $\mathrm{Red}_\Sigma(w) = \varnothing$.
\item The projective objects of $\mathrm{Sh}(C,G_\Sigma)$ are exactly the
coproducts of the $a\,y(u)$, $u \in \Sigma$.
\end{enumerate}
\end{corollary}

\begin{proof}
(a) $\mathrm{Sh}(C,G_\Sigma)$ is a presheaf topos by Theorem~\ref{res:6-6}, and
restriction carries $a\,y(u)$ to $y_{Q_\Sigma}(u)$ for $u \in \Sigma$:
indeed for a sheaf $F$, $\mathrm{Hom}(a\,y(u),F) \cong \mathrm{Hom}(y(u),
F) \cong F(u) \cong \mathrm{Hom}(y_{Q_\Sigma}(u), F|_\Sigma)$ by Yoneda on
both sides. Representables are projective and generate.

(b) By adjunction $\mathrm{Hom}(a\,y(w), F) \cong F(w)$ for every sheaf
$F$, and Proposition~\ref{res:6-5} identifies $F(w)$ with $\prod_{\pi}
F(s(\pi)) \cong \prod_\pi \mathrm{Hom}(a\,y(s(\pi)), F)$. A product of
corepresentable functors is corepresented by the coproduct of the
corepresenting objects, giving the claim; the case
$\mathrm{Red}_\Sigma(w) = \varnothing$ gives the empty coproduct.

(c) In any presheaf topos, projectives are retracts of coproducts of
representables. When the base category is idempotent-split these are
precisely the coproducts of representables; here the only idempotents of
$C|_\Sigma$ are identities (the quiver $Q_\Sigma$ is acyclic), which split
trivially.
\end{proof}

Part (b) is worth emphasising, since it is special to \emph{free}
categories: the splitting of matching families by final arrow in
Proposition~\ref{res:6-5} turns every sheafified representable into a coproduct.
Over a general site, $a\,y(w)$ need not be projective at all.

\subsection{The linear picture}\label{sub:6-6}

Write $A = kQ$, let $e = \sum_{u \in \Sigma} e_u$ be the corresponding
idempotent, and $B = eAe$.

\begin{proposition}\label{prop:linear}\label{res:6-8}
$B \cong kQ_\Sigma$. Moreover, let $\mathrm{Sh}_k \subseteq
\mathrm{mod}\text{-}A$ be the full subcategory of representations $M$ with
$M(w) \xrightarrow{\ \sim\ } \bigoplus_{a : v \to w} M(v)$ for every $w
\notin \Sigma$. Then restriction gives an equivalence $\mathrm{Sh}_k
\simeq \mathrm{mod}\text{-}kQ_\Sigma$, realized by the coinduction functor
\[
i_* = \mathrm{Hom}_B(Ae, -) : \mathrm{mod}\text{-}B \longrightarrow
\mathrm{mod}\text{-}A,
\]
the right adjoint of $i^* = (-)e$. Moreover $i^*$ is exact with kernel the
representations supported on $Q_0 \setminus \Sigma$, so
\[
\mathrm{mod}\text{-}kQ_\Sigma \;\simeq\; \mathrm{mod}\text{-}kQ \,/\,
\mathrm{mod}_{Q_0 \setminus \Sigma},
\]
a Serre quotient.
\end{proposition}

\begin{proof}
$B = eAe$ has $k$-basis the paths of $Q$ between $\Sigma$-vertices, and
$kQ_\Sigma$ has basis the words in arrows of $Q_\Sigma$; Lemma~\ref{res:6-4}'s
unique factorization is exactly a bijection between these bases,
compatible with concatenation. The defining condition of $\mathrm{Sh}_k$ is the one-step
form, Proposition~\ref{res:6-5}(iii), so $\mathrm{Sh}_k$ is the $k$-linear
sheaf category; the equivalence and its explicit inverse are then the
$k$-linear form of Theorem~\ref{res:6-6} and
Proposition~\ref{res:6-5}(ii), whose proofs use only limits (products) and
so apply verbatim over $k$; all products are finite because $Q$ is. That the right adjoint --- not the left ---
realizes the embedding is forced by the defining condition being a limit;
one sees it concretely on $Q = (1 \to 2)$ with $\Sigma = \{2\}$, where the
condition reads $M(1) = 0$, matching $\mathrm{Hom}_B(Ae_2, N) = (0,N)$ and
not $N \otimes_B e_2A$. Exactness of $(-)e$ and the identification of its
kernel are immediate.
\end{proof}

\begin{corollary}\label{res:6-9}
The indecomposable projectives of $\mathrm{Sh}_k \simeq
\mathrm{mod}\text{-}kQ_\Sigma$ are $P_{Q_\Sigma}(u) = k\,C(-,u)|_\Sigma$
for $u \in \Sigma$, and every projective is a direct sum of these. Since
$Q_\Sigma$ is acyclic, $kQ_\Sigma$ is hereditary for \emph{every} subset
$\Sigma \subseteq Q_0$.
\end{corollary}

\begin{proof}
Krull--Schmidt in $\mathrm{mod}\text{-}kQ_\Sigma$, which applies because
$kQ_\Sigma$ is finite-dimensional, plus Corollary~\ref{res:6-7}(c).
Heredity is the standard fact for path algebras of finite acyclic quivers.
\end{proof}

Two cautions. First, $a\,y(u)$ for $u \in \Sigma$ is \emph{not} in general
the projective $P_Q(u)$ of $\mathrm{mod}\text{-}kQ$ carried along: for $Q
= (1 \to 2)$ and $\Sigma = \{1\}$ the sheaf condition forces $M(2) \cong
M(1)$, so $a P_Q(1)$ has dimension vector $(1,1)$, which is $P_Q(2)$, not
$P_Q(1) = (1,0)$. The generators are \emph{indexed} by $\Sigma$, not given
by the old projectives at $\Sigma$. Second, $i_*$ is only left exact, so
$\mathrm{Ext}$ groups computed in $\mathrm{Sh}_k$ do not agree with those
computed in $\mathrm{mod}\text{-}kQ$; only $i^*$ is exact.

\subsection{Examples}\label{sub:6-7}

Take $Q$ equioriented of type $A_3$, $1 \xrightarrow{a} 2 \xrightarrow{b}
3$, as in \S\ref{sub:5-3}. The eight topologies give:

\[
\begin{array}{lll}
\Sigma & Q_\Sigma & \mathrm{Sh}(C,G_\Sigma) \\[2pt]
\hline
\{1,2,3\} & A_3 & \mathrm{PSh}(A_3) \\
\{1,2\} & A_2 & \mathrm{PSh}(A_2) \\
\{2,3\} & A_2 & \mathrm{PSh}(A_2) \\
\{1,3\} & A_2 \ (\text{via } ba) & \mathrm{PSh}(A_2) \\
\{1\},\ \{2\},\ \{3\} & \bullet & \mathbf{Set} \\
\varnothing & \text{empty} & \mathbf{1}
\end{array}
\]

For $\Sigma = \{1,3\}$ the reduced quiver has a single arrow, the
composite $ba$, whose interior vertex $2$ has been deleted --- this is the
sense in which $Q_\Sigma$ is a \emph{derived} quiver rather than a
subquiver. For $\Sigma = \{2\}$ one computes $J(1) = \varnothing$ (so
$F(1) = *$) and $J(3) = \{b, ba\}$ (so $F(3) \cong F(2)$), leaving
$\mathbf{Set}$, with sole projective generator $a\,y(2)$.

For the Kronecker quiver $1 \rightrightarrows 2$ of \S\ref{sub:5-2} with $\Sigma =
\{1\}$, we get $\mathrm{Red}_\Sigma(2) = \{a,b\}$, so the sheaf condition
reads $F(2) \cong F(1) \times F(1)$; the topos is $\mathbf{Set}$ and $a\,
y(2) \cong a\,y(1) \sqcup a\,y(1)$, illustrating Corollary~\ref{res:6-7}(b).

Finally, note the count: Corollary~\ref{res:4-3}'s $2^m$ topologies index precisely
the $2^m$ reduced quivers $Q_\Sigma$, i.e.\ the $2^m$ hereditary Serre
quotients of $\mathrm{mod}\text{-}kQ$ obtained by idempotent truncation.
Distinct subsets give distinct subcategories, though not always
inequivalent ones --- as the three singletons above show.

\section{Which subcategories arise}\label{sec:which}

Throughout $A = kQ$ with
$Q$ finite acyclic, and we work in $\mathrm{mod}\text{-}A$;
$\mathcal{S}_\Sigma \subseteq \mathrm{mod}\text{-}A$ denotes the
$k$-linear sheaf subcategory of Proposition~\ref{prop:linear}, so
$\mathcal{S}_\Sigma \simeq \mathrm{mod}\text{-}kQ_\Sigma$.

Theorem~\ref{thm:perp} is the specialization to this setting of
\cite[Thm.~4.2]{DLL}, which characterizes sheaves of modules on any
ringed site as the modules right perpendicular to the torsion modules; it
also follows from \cite[Prop.~5.1, Cor.~1.5(a)]{GL}, as detailed in
\S\ref{sec:literature};
we include the short direct proof because the projective resolution it
runs on is also what makes Lemma~\ref{lem:res} available for the rest of
the section. Propositions~\ref{prop:wide}--\ref{prop:lattice} and
\S\ref{sec:capture} record the position of $\mathcal{S}_\Sigma$ inside
the lattice of wide subcategories. These too are known, or immediate from
known results; \S\ref{sec:literature} gives the attributions, and the
proofs here are included because they are short and self-contained in the
quiver case.

\subsection{Perpendicular description}\label{sub:7-1}

Recall that for a class $\mathcal{X}$ of modules, the right perpendicular
category is
\[
\mathcal{X}^{\perp} = \{\, M \;\mid\; \mathrm{Hom}(X,M) = 0 =
\mathrm{Ext}^1(X,M) \ \text{ for all } X \in \mathcal{X} \,\}.
\]

\begin{lemma}\label{lem:res}\label{res:7-1}
For every vertex $v$ there is an exact sequence
\[
0 \longrightarrow \bigoplus_{a\,:\,u \to v} P(u) \longrightarrow P(v)
\longrightarrow S_v \longrightarrow 0,
\]
and consequently, for every $M$,
\[
\mathrm{Hom}(S_v, M) = \ker \theta^M_v, \qquad \mathrm{Ext}^1(S_v,M) =
\operatorname{coker} \theta^M_v,
\]
where $\theta^M_v : M(v) \to \bigoplus_{a : u \to v} M(u)$ is the map
assembled from the structure maps of $M$.
\end{lemma}

\begin{proof}
$P(v) = k\,C(-,v)$ has basis the paths into $v$, and the kernel of $P(v)
\to S_v$ is spanned by the paths of length $\geq 1$ into $v$. Sorting
these by final arrow gives a decomposition $\bigoplus_{a : u \to v}
a\cdot P(u)$, and $g \mapsto ag$ is a bijection from paths into $u$ onto
paths into $v$ with final arrow $a$ --- injective and with the stated
image because $C$ is a free path category. This is the displayed
sequence. Applying $\mathrm{Hom}(-,M)$, using $\mathrm{Hom}(P(x),M) =
M(x)$ (Yoneda) and $\mathrm{Ext}^{\geq 1}(P(x),M) = 0$, yields
\[
0 \to \mathrm{Hom}(S_v,M) \to M(v) \xrightarrow{\ \theta^M_v\ }
\bigoplus_{a : u \to v} M(u) \to \mathrm{Ext}^1(S_v,M) \to 0,
\]
which is the claim.
\end{proof}

\begin{theorem}\label{thm:perp}\label{res:7-2}
For every $\Sigma \subseteq Q_0$,
\[
\mathcal{S}_\Sigma \;=\; \Big( \bigoplus_{v \notin \Sigma} S_v
\Big)^{\perp} \;=\; \bigcap_{v \notin \Sigma} S_v^{\perp}.
\]
Thus the Grothendieck topologies of Corollary~\ref{res:4-3} correspond exactly to
the perpendicular categories of sets of simple modules.
\end{theorem}

\begin{proof}
By Lemma~\ref{res:7-1}, $M \in S_v^{\perp}$ iff $\theta^M_v$ is an isomorphism. By
Proposition~\ref{res:6-8}, $M \in \mathcal{S}_\Sigma$ iff $\theta^M_w$ is an
isomorphism for every $w \notin \Sigma$.
\end{proof}

Since $Q$ is acyclic, any linear extension of the reachability order makes
$\{S_v\}_{v \notin \Sigma}$ into an exceptional sequence, so Theorem~\ref{res:7-2}
places $\mathcal{S}_\Sigma$ inside the Geigle--Lenzing theory of
perpendicular categories \cite{GL}; in particular each $\mathcal{S}_\Sigma$ is
functorially finite in $\mathrm{mod}\text{-}A$, and $D^b(\mathcal{S}_\Sigma)$
embeds as a thick subcategory of $D^b(\mathrm{mod}\text{-}A)$. We do not
use these facts below.

\subsection{They are wide, and rarely Serre}\label{sub:7-2}

\begin{proposition}\label{prop:wide}\label{res:7-3}
Each $\mathcal{S}_\Sigma$ is a wide subcategory of $\mathrm{mod}\text{-}A$:
it is closed under direct sums and summands, and under kernels, cokernels
and extensions taken in $\mathrm{mod}\text{-}A$. Its inclusion is exact,
and $\mathcal{S}_\Sigma \simeq \mathrm{mod}\text{-}kQ_\Sigma$ is again
hereditary.
\end{proposition}

\begin{proof}
Fix $w \notin \Sigma$. Both $F : M \mapsto M(w)$ and $G : M \mapsto
\bigoplus_{a : u \to w} M(u)$ are exact functors $\mathrm{mod}\text{-}A
\to \mathrm{mod}\text{-}k$, and $\theta_w : F \Rightarrow G$ is a natural
transformation; $\mathcal{S}_\Sigma$ is the full subcategory on which
$\theta_w$ is invertible for all such $w$.

If $f : M \to N$ has $\theta_w$ invertible at $M$ and $N$, then exactness
of $F$ and $G$ identifies $F(\ker f) = \ker F(f)$ and $G(\ker f) = \ker
G(f)$, and naturality identifies $\theta_w$ at $\ker f$ with the
restriction of an isomorphism between these kernels; the same argument
with cokernels handles $\operatorname{coker} f$. For an extension $0 \to M
\to X \to N \to 0$, applying $F$ and $G$ gives short exact sequences and
the five lemma makes $\theta_w$ at $X$ an isomorphism. Direct sums and
summands are clear. Heredity is Corollary~\ref{res:6-9}.
\end{proof}

\begin{proposition}\label{prop:serre}\label{res:7-4}
$S_v \in \mathcal{S}_\Sigma$ if and only if $v \in \Sigma$ and no path in
$Q$ leads from $v$ to a vertex outside $\Sigma$. Consequently
$\mathcal{S}_\Sigma$ is a Serre subcategory of $\mathrm{mod}\text{-}A$ if
and only if $\Sigma$ is closed under successors, in which case $Q_\Sigma$
is the full subquiver on $\Sigma$ and $\mathcal{S}_\Sigma$ is precisely
the category of representations supported on $\Sigma$.
\end{proposition}

\begin{proof}
For the first claim, test the condition of Proposition~\ref{res:6-5}. If $v \notin
\Sigma$, take $w = v$: the left side is $S_v(v) = k$ while every $\pi \in
\mathrm{Red}_\Sigma(v)$ has $s(\pi) \in \Sigma$, so $S_v(s(\pi)) = 0$ and
the right side vanishes; hence $S_v \notin \mathcal{S}_\Sigma$. If $v \in
\Sigma$, then for $w \notin \Sigma$ the left side is $0$ and the right
side is $\bigoplus_{\pi : s(\pi) = v} k$, so the condition holds for all
$w \notin \Sigma$ exactly when no $\Sigma$-reduced path runs from $v$ to a
vertex outside $\Sigma$ --- equivalently, when no path at all does, since
a path from $v$ leaving $\Sigma$ has an initial segment that is
$\Sigma$-reduced and ends outside $\Sigma$.

Now suppose $\mathcal{S}_\Sigma$ is Serre. Serre subcategories of the
length category $\mathrm{mod}\text{-}A$ are exactly those generated by a
set of simples \cite[Prop.~5.3]{GL}, so $\mathcal{S}_\Sigma$ is generated by $T = \{v \mid S_v
\in \mathcal{S}_\Sigma\}$ and therefore has exactly $|T|$ simple objects.
But $\mathcal{S}_\Sigma \simeq \mathrm{mod}\text{-}kQ_\Sigma$ has
$|\Sigma|$ simple objects, so $|T| = |\Sigma|$, forcing $T = \Sigma$ by
the first claim --- i.e.\ $\Sigma$ is successor-closed. Conversely if
$\Sigma$ is successor-closed then $\mathrm{Red}_\Sigma(w) = \varnothing$
for every $w \notin \Sigma$, so by Proposition~\ref{res:6-5} $\mathcal{S}_\Sigma =
\{M \mid M(w) = 0 \text{ for } w \notin \Sigma\}$, which is Serre, and no
$\Sigma$-reduced path has interior vertices, so $Q_\Sigma$ is the full
subquiver.
\end{proof}

So the generic member of the family is wide but \emph{not} Serre, and is
also neither a torsion class nor a torsion-free class: for $Q = (1 \to 2)$
and $\Sigma = \{1\}$ the subcategory is $\{\,\phi : V_2 \to V_1 \text{
invertible}\,\}$, which contains $(k \xrightarrow{\ \sim\ } k)$ but
neither its subobject $(0 \to k)$ nor its quotient $(k \to 0)$.

\subsection{The lattice}\label{sub:7-3}

\begin{proposition}\label{prop:lattice}\label{res:7-5}
The assignment $\Sigma \mapsto G_\Sigma$ is an anti-isomorphism from the
Boolean lattice $2^{Q_0}$ onto the lattice of Grothendieck topologies on
$C$, and $\Sigma \mapsto \mathcal{S}_\Sigma$ is an injective
order-preserving map into the lattice of wide subcategories of
$\mathrm{mod}\text{-}A$.
\end{proposition}

\begin{proof}
Topologies are ordered by $G \leq G'$ iff $G(v) \subseteq G'(v)$ for all
$v$, i.e.\ iff $J_{G'}(v) \subseteq J_G(v)$ for all $v$. If $\Sigma
\subseteq \Sigma'$ and $v \in \Sigma$ then $J_{\Sigma}(v) = t_v =
J_{\Sigma'}(v)$; if $v \notin \Sigma'$ then Lemma~\ref{res:6-1} gives $J_\Sigma(v)
\subseteq J_{\Sigma'}(v)$ since meeting $\Sigma$ implies meeting
$\Sigma'$; and if $v \in \Sigma' \setminus \Sigma$ then $J_{\Sigma'}(v) =
t_v \supseteq J_\Sigma(v)$. So $\Sigma \subseteq \Sigma'$ implies
$G_{\Sigma'} \leq G_\Sigma$. Conversely if $G_{\Sigma'} \leq G_\Sigma$ and
$v \in \Sigma$, then $t_v = J_\Sigma(v) \subseteq J_{\Sigma'}(v)$ forces
$J_{\Sigma'}(v) = t_v$, i.e.\ $v \in \Sigma'$. With Corollary~\ref{res:4-3} this is
an anti-isomorphism.

Monotonicity of $\Sigma \mapsto \mathcal{S}_\Sigma$ is Theorem~\ref{res:7-2}:
enlarging $\Sigma$ deletes conditions. For injectivity, note
$\mathcal{S}_\Sigma$ is reflective in $\mathrm{mod}\text{-}A$ with
reflector $i_*i^*$, and a reflector is determined up to isomorphism by the
subcategory; hence so is the class $\{M \mid i_*i^*M = 0\} = \{M \mid Me =
0\}$, the modules supported on $Q_0 \setminus \Sigma$, which determines
$\Sigma$.
\end{proof}

That the lattice of topologies is \emph{Boolean} is worth a remark: on a
general site the topologies form a frame with no reason to be
complemented. Here complementation $\Sigma \mapsto Q_0 \setminus \Sigma$
is visible on the nose.

\subsection{How much of the wide lattice is captured}\label{sec:capture}

Not all wide subcategories arise. By Ingalls--Thomas \cite{IT}, for $Q$
of Dynkin type $A_n$ the wide subcategories of $\mathrm{mod}\text{-}kQ$
are in bijection with the noncrossing partitions of $S_{n+1}$, hence are
counted by the Catalan number $C_{n+1}$, whereas our family has $2^n$
members.

The smallest case is already instructive. For $Q = (1 \to 2)$ the
indecomposables are $S_1 = P(1)$, $S_2$, and $P = P(2)$, and there are
$C_3 = 5$ wide subcategories. The four arising from topologies are:
\[
\begin{array}{lll}
\Sigma & \mathcal{S}_\Sigma & \text{type} \\[2pt]
\hline
\{1,2\} & \mathrm{mod}\text{-}kQ & \text{everything} \\
\{2\} & \langle S_2 \rangle = \{M(1) = 0\} & \text{Serre} \\
\{1\} & \mathrm{add}(P) = \{\theta \text{ invertible}\} & \text{wide, not Serre} \\
\varnothing & 0 & \text{Serre}
\end{array}
\]
The missing fifth is $\langle S_1 \rangle$, the modules supported at
vertex $1$. It is Serre, but $\{1\}$ is not successor-closed, and the
subset $\{1\}$ is already spent on $\mathrm{add}(P)$. This is the general
pattern: our family sees each subset once, and by
Proposition~\ref{prop:serre} it meets the Serre subcategories exactly in
the successor-closed subsets --- for the equioriented $A_n$ that is only
$n+1$ of the $2^n$ members.

\section{Relation to the literature}\label{sec:literature}

Because the classification of topologies used here has been proved
repeatedly and at increasing generality, we set out carefully what is due
to whom.

\subsection{The rigidity ladder}\label{sub:8-1}

A Grothendieck topology is \emph{rigid}, or a \emph{subcategory topology},
when every object is covered by the sieve generated by morphisms to
irreducible objects; equivalently, when the associated subtopos of
$\mathrm{PSh}(C)$ is again a presheaf topos, on the full subcategory of
irreducibles. The following is the sequence of known results, from most
to least general.

\begin{itemize}
\item Grothendieck and Verdier \cite[IV.9]{SGA4}: every Grothendieck
topology on a finite Karoubi category is rigid. See also
\cite[C2.2.21]{Elephant}, and \cite{KL} for the lattice of essential
localizations.
\item Di, Li and Liang \cite[Thm.~5.7]{DLL}: for a directed category
$\mathcal{C}$ with every $\mathcal{C}(x,x)$ finite, every Grothendieck
topology on $\mathcal{C}^{\mathrm{op}}$ is rigid if and only if
$\mathcal{C}$ is a noetherian EI category.
\item Lindenhovius \cite{Lindenhovius} for artinian posets; Hemelaer
\cite{Hemelaer} for general posets, via locale and domain theory.
\item Murfet \cite{Murfet} for finite acyclic quivers.
\end{itemize}

The path category of a finite acyclic quiver is a finite category whose
only endomorphisms are identities, so it satisfies the hypotheses of all
four simultaneously.

\subsection{A dictionary of conventions}

Comparison with \cite{DLL} requires care with variance, so we record the
dictionary. Our sieves consist of morphisms \emph{into} an object, so a
Grothendieck topology in our sense lives on $C$ and its sheaves are
presheaves on $C$. In \cite{DLL} a sieve on $x$ is a left ideal of
morphisms \emph{out of} $x$ --- a sieve on $x$ in $\mathcal{C}^{\mathrm{op}}$
--- and representations are covariant functors $\mathcal{C} \to
\mathbf{Set}$. The two conventions match under
\[
\mathcal{C} = C^{\mathrm{op}},
\]
so that their $\mathcal{C}^{\mathrm{op}}$ is our $C$, their sieves are
ours, and their representations of $\mathcal{C}$ are our presheaves on
$C$. Since $C$ is the path category of $Q$, their $\mathcal{C}$ is the
path category of $Q^{\mathrm{op}}$, which is again finite and acyclic;
all the cited statements are invariant under this relabelling. The same
remark applies to \cite{Hemelaer} and \cite{Lindenhovius}.

\subsection{What is not ours}\label{sub:8-2}

\begin{itemize}
\item \S\S\ref{sec:quivers}--\ref{sec:bijection} are Murfet's
\cite{Murfet}. His consistency condition (our \eqref{eq:consistent}) appears
verbatim as condition (5.1) of \cite{DLL}, where the authors state that
they are generalizing his construction. Parts (a) and (b) of our
Lemma~\ref{res:3-2} are \cite[Lem.~2.2]{DLL}, and its consequence that a
minimal covering sieve exists corresponds to \cite[Lem.~5.15]{DLL}, which
proves closure under \emph{arbitrary} intersections in the noetherian EI
setting; our Theorem~\ref{res:4-1} is \cite[Prop.~5.16]{DLL}.
\item Theorem~\ref{thm:comparison} is the Comparison Lemma
\cite[C2.2.3]{Elephant} applied to a rigid topology; in the form used
here it is \cite[Thm.~5.18]{DLL} and \cite[Cor.~1.6]{DLL}.
\item Theorem~\ref{thm:perp} is the specialization of
\cite[Thm.~4.2]{DLL}, which shows for an arbitrary ringed site that a
presheaf of modules is a sheaf if and only if it is right perpendicular,
in the sense of Geigle--Lenzing \cite{GL}, to the torsion presheaves.
The Serre quotient description in Proposition~\ref{prop:linear} is
\cite[Cor.~1.2]{DLL}.
\item The example of \S\ref{sec:examples} on the Kronecker quiver
appears, with the same four topologies, as \cite[Ex.~3.16]{DLL}.
\item The correspondence between topologies and hereditary torsion pairs
is the category-algebra analogue of the classical correspondence between
Gabriel topologies on a ring and hereditary torsion theories on its
module category \cite{Stenstrom}; this is noted in \cite[\S3]{DLL}.
\end{itemize}

\subsection{Attributions for \S\ref{sec:which}}\label{sub:8-3}

The results of \S\ref{sec:which} are likewise known. We give the proofs
because in the quiver case they are short and entirely combinatorial, but
we claim no priority.

\begin{itemize}
\item Proposition~\ref{prop:wide} (that each $\mathcal{S}_\Sigma$ is
wide) is \cite[Prop.~1.1]{GL}: for a class $\mathcal{X}$ of objects of an
abelian category, $\mathcal{X}^{\perp}$ is always closed under kernels
and extensions, and is an exact abelian subcategory whenever
$\operatorname{proj\,dim} \mathcal{X} \leq 1$ --- which holds for our
$\mathcal{X} = \{S_v\}_{v \notin \Sigma}$ since $kQ$ is hereditary. Our
proof, via invertibility of a natural transformation between exact
functors, is a rediscovery of that statement in this special case.

\item Theorem~\ref{thm:perp}, Corollary~\ref{res:6-9} and the Serre
quotient description are all contained in \cite[\S5]{GL}. Take $P = eA$
with $e = \sum_{u \in \Sigma} e_u$, so $\mathrm{End}(P) = eAe$. Then
$P^{\perp}$ is the Serre subcategory of modules annihilated by the trace
ideal of $P$ --- the modules supported off $\Sigma$ --- and
\cite[Prop.~5.1]{GL} gives
\[
(P^{\perp})^{\perp} \simeq \mathrm{mod}\text{-}eAe
\qquad\text{and}\qquad
\mathrm{mod}\text{-}A / P^{\perp} \simeq \mathrm{mod}\text{-}eAe;
\]
by \cite[Cor.~1.5(a)]{GL} the perpendicular of the Serre subcategory
$P^{\perp}$ equals the perpendicular of its simples, i.e.\
$(P^{\perp})^{\perp} = \big(\bigoplus_{v\notin\Sigma} S_v\big)^{\perp} =
\mathcal{S}_\Sigma$. The functor-category form \cite[Prop.~5.1$^*$]{GL}
is the $k$-linear comparison theorem; Geigle--Lenzing explicitly
identify these two propositions with the representation-theoretic
\emph{deletion of vertices}.

\item Proposition~\ref{prop:serre} rests on \cite[Prop.~5.3]{GL}, which
proves for any Artin algebra that the Serre subcategories of
$\mathrm{mod}$ are exactly those generated by a set of simples,
equivalently the perpendicular categories $P^{\perp}$ of projectives,
equivalently the subcategories $\mathrm{mod}(\Lambda/\mathfrak{a})$ for
idempotent ideals $\mathfrak{a}$. Our successor-closed criterion is the
quiver-combinatorial translation. (For the $\tau$-perpendicular
generalization of this circle of ideas see \cite{Jasso}.)

\item The count $2^{|Q_0|}$ has an algebra-side ancestor:
\cite[Cor.~5.4]{GL} shows that an Artin algebra with $n$ simple modules
admits at most $2^n$ surjective homological epimorphisms, because
$\mathrm{mod}$ has exactly $2^n$ Serre subcategories.

\item Ingalls and Thomas \cite{IT} showed that for $Q$ acyclic the wide
subcategories of $\mathrm{mod}\text{-}kQ$ are in bijection with
noncrossing partitions, clusters and functorially finite torsion classes.
Proposition~\ref{prop:lattice} and \S\ref{sec:capture} locate the sheaf
subcategories inside that classification; the count is an immediate
consequence of Corollary~\ref{res:4-3} and their theorem.
\end{itemize}

Two remarks on scope. First, Proposition~\ref{prop:lattice}'s Boolean
lattice is really a statement about idempotents rather than about quivers:
for a general rigid topology the irreducibles form a retract-closed full
subcategory, and topologies correspond to retract-closed subcategories;
when the only idempotents are identities, as in an acyclic path category,
this collapses to subsets. Second, Proposition~\ref{prop:serre}'s
criterion is stated in terms of successor-closed sets of vertices, which
in the language of \cite[Def.~5.3]{DLL} is the condition that $\Sigma$ be
an ideal of the directed category; we expect the criterion to hold
verbatim for noetherian EI categories, with \cite[Thm.~5.18(1)]{DLL} as
the input, but we have not checked this.

\subsection{The intrinsic characterization}\label{sub:8-4}

This is the intrinsic form of Question~\ref{res:1-1}.
Theorem~\ref{thm:perp} identifies the sheaf subcategories as the
perpendicular categories of sets of \emph{simples}, a proper subfamily of
the perpendicular categories of arbitrary exceptional sequences: for $Q$
of type $A_2$ the subfamily has $4$ of the $5$ wide subcategories
(\S\ref{sec:capture}), and in general it has $2^{|Q_0|}$ members. The
machinery of \cite{GL} recognizes the members of the subfamily
intrinsically.

\begin{proposition}\label{prop:intrinsic}
A full subcategory $\mathcal{W} \subseteq \mathrm{mod}\text{-}kQ$ is a
sheaf subcategory --- i.e.\ $\mathcal{W} = \mathcal{S}_\Sigma$ for some
$\Sigma \subseteq Q_0$ --- if and only if
\[
\mathcal{W} = ({}^{\perp}\mathcal{W})^{\perp}
\qquad\text{and}\qquad
{}^{\perp}\mathcal{W} \text{ is a Serre subcategory.}
\]
\end{proposition}

\begin{proof}
($\Rightarrow$) Let $\mathcal{W} = \mathcal{S}_\Sigma =
\mathcal{S}^{\perp}$, where $\mathcal{S}$ is the Serre subcategory of
modules supported off $\Sigma$ (Theorem~\ref{thm:perp} together with
\cite[Cor.~1.5(a)]{GL}). Serre subcategories of $\mathrm{mod}$ of an
Artin algebra are localizing \cite[Prop.~5.3, (i)$\Leftrightarrow$(iii)]{GL}, and for a
localizing subcategory $\mathcal{S} = {}^{\perp}(\mathcal{S}^{\perp})$
\cite[Cor.~2.3]{GL}. Hence ${}^{\perp}\mathcal{W} = \mathcal{S}$ is
Serre, and $({}^{\perp}\mathcal{W})^{\perp} = \mathcal{S}^{\perp} =
\mathcal{W}$.

($\Leftarrow$) Suppose ${}^{\perp}\mathcal{W} =: \mathcal{S}$ is Serre.
Then $\mathcal{S}$ is generated by the simples it contains
\cite[Prop.~5.3, (i)$\Leftrightarrow$(ii)]{GL}, say $\{S_v\}_{v \notin \Sigma}$ for some subset
$\Sigma$, and by \cite[Cor.~1.5(a)]{GL} again, $\mathcal{W} =
\mathcal{S}^{\perp} = \big(\bigoplus_{v \notin \Sigma}
S_v\big)^{\perp} = \mathcal{S}_\Sigma$.
\end{proof}

For the missing fifth wide subcategory of \S\ref{sec:capture} this test
is concrete: $\mathcal{W} = \langle S_1 \rangle$ has
${}^{\perp}\mathcal{W} = \mathrm{add}(P)$, which is not Serre, so
$\mathcal{W}$ is not a sheaf subcategory --- while for $\mathcal{W} =
\mathrm{add}(P)$ one finds ${}^{\perp}\mathcal{W} = \langle S_2
\rangle$, which is Serre. Note the asymmetry: the test uses the
\emph{left} perpendicular, and $\langle S_1 \rangle$ fails it even
though $\langle S_1 \rangle$ is itself Serre --- being Serre is about
$\mathcal{W}$'s position as a kernel, being a sheaf subcategory about its
position as a perpendicular.

\end{document}